\documentclass[12pt,english]{smfart}
\usepackage{smfenum}
\usepackage{bull}
\usepackage{amscd}
\usepackage{amssymb}
\usepackage{epsfig}
\usepackage{graphicx}
\usepackage[all]{xy}

\advance\headheight 2pt
\calclayout

\usepackage[T1]{fontenc}
\usepackage{ae}

\renewcommand\AA{\mathbb{A}}
\newcommand\QQ{\mathbb{Q}}

\newcommand\NN{\mathcal{N}}
\newcommand\RR{\mathbb{R}}
\newcommand\R{\mathcal{R}}
\newcommand\TT{\mathcal{T}}
\newcommand\ZZ{\mathbb{Z}}
\newcommand\ZZp{\ZZ_{>0}}
\newcommand\PP{\mathbb{P}}
\newcommand\xx{\mathbf{x}}
\newcommand\dd{\,\mathrm{d}}
\newcommand\Pone{{\PP^1}}
\newcommand\Ptwo{{\PP^2}}
\newcommand\Pthree{{\PP^3}}
\newcommand\Pfour{{\PP^4}}
\newcommand{\sbase}[6]{\eta^{({#1},{#2},{#3},{#4},{#5},{#6})}}
\newcommand{\congr}[3]{{#1} \equiv {#2}\ (\mathrm{mod}\ {#3})}
\newcommand\ep{\epsilon}
\newcommand\lle{{\ll_\ep}}
\newcommand\NU{{N(\So,-K_S,B)}}
\newcommand\Res{\mathrm{Res}}
\newcommand\prim{\mathrm{prim}}
\newcommand{\Aone}{{\mathbf A}_1}
\newcommand{\Atwo}{{\mathbf A}_2}
\newcommand{\Athree}{{\mathbf A}_3}
\newcommand{\Afour}{{\mathbf A}_4}
\newcommand{\Afive}{{\mathbf A}_5}
\newcommand{\Dfour}{{\mathbf D}_4}
\newcommand{\Dfive}{{\mathbf D}_5}
\newcommand{\Esix}{{\mathbf E}_6}
\newcommand\aone{\alpha_1}
\newcommand\atwo{\alpha_2}
\newcommand\athr{\alpha_3}
\newcommand{\Fone}{\eta_2}
\newcommand{\Ftwo}{\eta_3\eta_5^2}
\newcommand{\Fthr}{\eta_4\eta_6^2}
\newcommand{\Tone}{\aone^2\Fone}
\newcommand{\Ttwo}{\atwo\Ftwo}
\newcommand{\Tthr}{\athr\Fthr}
\newcommand{\Cone}{\eta_1\eta_3\eta_4\eta_5\eta_6}
\newcommand{\Ctwo}{\eta_1\eta_2\eta_3\eta_4\eta_6}
\newcommand{\Cthr}{\eta_1\eta_2\eta_3\eta_4\eta_5}
\newcommand{\tS}{{\widetilde S}}
\newcommand{\So}{{S^\circ}}
\newcommand{\anti}{-K_\tS}
\renewcommand{\le}{\leqslant}
\renewcommand{\ge}{\geqslant}
\newcommand{\ee}{\boldsymbol{\eta}}
\renewcommand{\aa}{\boldsymbol{\alpha}}
\newcommand{\ai}{\alpha_i}
\newcommand{\nd}{\nmid}
\newcommand{\restr}{\text{restriction}}
\newcommand{\autom}{\text{automatically}}
\newcommand{\allow}{\text{allowed}}
\newcommand{\qa}[1]{q_{#1,1}}
\newcommand{\qb}[1]{q_{#1,2}}
\newcommand{\qc}[1]{q_{#1,3}}

\DeclareMathOperator{\NS}{NS}
\DeclareMathOperator{\Pic}{Pic}
\DeclareMathOperator{\Spec}{Spec}
\DeclareMathOperator{\rk}{rk}
\DeclareMathOperator{\additiv}{a}
\newcommand\Ga{\mathbb{G}_{\additiv}}
\DeclareMathOperator{\multiplikativ}{m}
\newcommand\Gm{\mathbb{G}_{\multiplikativ}}

\newcommand\rto{\dashrightarrow}

\newtheorem{theorem}{Theorem}
\newtheorem{lemma}[theorem]{Lemma}

\newtheorem{conj}[theorem]{Conjecture}
\theoremstyle{definition}

\newtheorem{example}[theorem]{Example}

\numberwithin{equation}{section}

\begin{document}

\title[Torsors and rational points]{Universal torsors over Del Pezzo surfaces and rational points}
\author{Ulrich Derenthal}
\address{Mathematisches Institut\\Universit\"at G\"ottingen\\
  Bunsenstr. 3-5\\37073 G\"ottingen\\Germany}
\email{derentha@math.uni-goettingen.de}
\author{Yuri Tschinkel}
\address{Mathematisches Institut\\Universit\"at G\"ottingen\\
  Bunsenstr. 3-5\\37073 G\"ottingen\\Germany}
\curraddr{Courant Institute of Mathematical Sciences\\New York University\\
  251 Mercer St.\\New York, NY 10012}
\email{yuri@uni-math.gwdg.de\\tschinkel@cims.nyu.edu}
\date{April 8, 2006}

\begin{abstract}
  We discuss Manin's conjecture concerning the distribution of
  rational points of bounded height on Del Pezzo surfaces, and its
  refinement by Peyre, and explain applications of universal torsors to
  counting problems.  To illustrate the method, we provide a proof of
  Manin's conjecture for the unique split singular quartic Del Pezzo
  surface with a singularity of type $\Dfour$.
\end{abstract}

\maketitle

\tableofcontents

\section{Introduction}

Let $f\in \ZZ[x_0,\ldots, x_n]$ be a non-singular form of degree $d$. By the
circle method,
\[N(f,B):= \#\{ \xx \in \ZZ^{n+1}/\pm \mid \max_j(|x_j|)\le B\} \sim
c\cdot B^{n+1-d}\] (where $\xx \in \ZZ^{n+1}/\pm$ means that we
identify $\xx$ with $-\xx = (-x_0, \dots, -x_n)$) with $c\in
\RR_{>0}$, provided that $n \ge 2^d\cdot(d-1)$, and $f(\xx) = 0$ has
solutions over all completions of $\QQ$ (see \cite{MR0150129}).  Let
$X=X_f\subset \PP^n$ be a smooth hypersurface over $\QQ$, given by
$f(\xx)=0$.  It follows that
\begin{equation}\label{eqn:linear}
N(X,-K_X,B)\:=\#\{ \xx\in X(\QQ) \mid H_{-K_X}(\xx)\le B\} \sim C\cdot B,
\end{equation}
as $B\to \infty$. Here $X(\QQ)$ is the set of rational points on $X$,
represented by primitive vectors $\xx \in (\ZZ^{n+1}_\prim \setminus
0) / \pm$ (i.e., $\xx = (x_0, \dots, x_n)$ is identified with $-\xx$,
and there is no prime dividing all coordinates $x_0, \dots, x_n$), and
\begin{equation}\label{eqn:height}
H_{-K_X}(\xx):=
\max_j(|x_j|)^{n+1-d}, \,\, \text{ for }\,\, \xx=(x_0,\dots, x_n) \in
(\ZZ^{n+1}_\prim\setminus 0)/\pm.
\end{equation}
is the \emph{anticanonical height} of a primitive representative.

In 1989 Manin initiated a program towards understanding connections
between certain geometric invariants of algebraic varieties over
number fields and their arithmetic properties, in particular, the
distribution of rational points of bounded height, see
\cite{MR89m:11060} and \cite{MR1032922}.  The main goal is an
extension of the asymptotic formula \eqref{eqn:linear} to other
algebraic varieties of \emph{small} degree, called Fano varieties,
which are not necessarily isomorphic to hypersurfaces in projective
space.  

It became apparent, that in general, to obtain a geometric
interpretation of asymptotic results, it may be necessary to restrict
to appropriate Zariski open subsets of $X$. Otherwise, the number of
rational points on a Zariski closed subset of lower dimension may
dominate the total number of rational points; e.g., this phenomenon occurs for the surface \eqref{eq:surface} below where we will restrict
to the complement of its lines. Furthermore, we may need to allow
finite field extensions: while $X(\QQ)$ might be empty, $X(k)$ could
still contain infinitely many points for some number field $k$.

Of particular interest are Del Pezzo surfaces (cf. \cite{MR833513}),
e.g., cubic surfaces $S_3\subset \PP^3$ or degree 4 surfaces
$S_4:=Q_1\cap Q_2\subset \Pfour$, where $Q_1, Q_2$ are \emph{quadrics}
(defined by homogeneous equations of degree 2 in $x_0, \dots, x_4$).
Geometrically, smooth Del Pezzo surfaces are obtained by blowing up
$\le 8$ \emph{general points}\footnote{no three points on a line, no
six points on a curve of degree 2, no eight points with one of them
singular on a curve of degree 3} in $\Ptwo$. Blowing up is a standard
procedure in algebraic geometry (cf. \cite[Section
I.4]{MR0463157}). The blow-up $\pi: S' \to S$ of a surface $S$ at a
point $p$ replaces $p$ by a curve $E$ in a particular way. We have $S
\setminus\{p\} \cong S' \setminus E$, so $S$ and $S'$ are
\emph{birational}.  In our situation, this shows that Del Pezzo
surfaces are birational to $\Ptwo$, provided the ground field is
algebraically closed.

We can think of \emph{divisors} on blow-ups $S$ of $\Ptwo$ as formal
sums of curves on $S$. Considering divisors up to a certain
equivalence relation (see \cite[Section II.6]{MR0463157}) leads to the
\emph{Picard group} $\Pic(S)$ of divisor classes on $S$. 

For two curves on $S$ which intersect transversally, their
\emph{intersection number} is the number of intersection points. As
explained in \cite[Section V.1]{MR0463157}, this can be extended to
arbitrary divisor classes, defining the non-degenerate
\emph{intersection form} $(\cdot, \cdot)$ on $\Pic(S)$. In particular,
this defines the \emph{self intersection number} $(E,E)$ of (the class
of) a curve $E$. Of special interest are irreducible curves for which
this number is negative. We call them \emph{exceptional curves}. For
smooth Del Pezzo surfaces of degree 3 and 4, the exceptional curves
are exactly the lines (in the standard embedding considered above),
having self intersection number $-1$.

The singular Del Pezzo surfaces are obtained as follows: we blow up
$\Ptwo$ in special configurations of points (e.g., three points on a
line). This results in a smooth surface $\tS$ containing exceptional
curves with self intersection number $-2$ (called $(-2)$-curves; we do
not permit to blow up points on $(-2)$-curves subsequently).
Contracting the $(-2)$-curves gives a singular Del Pezzo surface $S$
whose \emph{minimal desingularization} is $\tS$. For the
surface~\eqref{eq:surface} below, more details can be found in
Section~\ref{sec:geometry}.

For number fields, we say that a Del Pezzo surface is split if all of
the exceptional curves are defined over that ground field;
there exist \emph{non-split} forms, some of which are not birational
to $\Ptwo$ over that ground field.

From now on, we work over $\QQ$. Manin's conjecture
in the special case of Del Pezzo surfaces can be formulated as follows.

\begin{conj}\label{conj:manin}
  Let $S$ be a Del Pezzo surface with at most rational double
  points\footnote{``mild'' singularities which can be resolved by
    blow-ups to a curve whose irreducible components are isomorphic to
    $\Pone$} over $\QQ$. Then there exists a subset $\So \subset S$
  which is dense and open in the Zariski topology such that
  \begin{equation}\label{eqn:asym}
    N(\So, -K_S,B)\sim c_{S,H}\cdot B(\log B)^{r-1},
  \end{equation}
  as $B\rightarrow \infty$, where $r$ is the rank of the
  Picard group of the minimal desingularization $\tS$ of $S$,
  over $\QQ$.
\end{conj}

The constant $c_{S,H}$ has been defined by Peyre \cite{MR1340296}; it
should be non-zero if $S(\QQ)\neq \emptyset$.  Note that a line
defined over $\QQ$ on a Del Pezzo surface such as $S_3$ or $S_4$
contributes $\sim c\cdot B^2$ rational points to the counting function
(for some positive constant $c$). Thus it is expected that $\So$ is
the complement to all lines defined over $\QQ$ (exceptional curves).

Table~\ref{tab:overview} gives an overview of current results towards
Conjecture~\ref{conj:manin} for Del Pezzo surfaces. In Column~4
(``type of result''), ``asymptotic'' means that the analog of
\eqref{eqn:asym} is established, including the predicted value of the
constant; ``bounds'' means that only upper and lower bounds of the
expected order of magnitude $B(\log B)^{r-1}$ with unknown constants
are proved.

The paper \cite{MR1620682} contains a proof of Manin's conjecture for
toric Fano varieties, including all smooth Del Pezzo surfaces of
degree $\ge 6$ and the $3\Atwo$ cubic surface\footnote{Singular
  Del Pezzo surfaces $S$ will be labeled by the type (in the
  ADE-classification) and number of their singularities. The
  corresponding Dynkin diagram describes the number and intersection
  behaviour of the $(-2)$-curves on $\tS$.}.  This result also covers:
\begin{itemize}
\item all singular surfaces of degree $\ge 7$ (i.e., $\Aone$ in degree
  7 and 8),
\item $\Aone$, $2\Aone$, $\Atwo+\Aone$ in degree 6,
\item $2\Aone$, $\Atwo+\Aone$ in degree 5,
\item $4\Aone$, $\Atwo+2\Aone$, $\Athree+2\Aone$ in degree 4.
\end{itemize}

\begin{table}[htbp]\centering
\begin{tabular}[ht]{|c|c|c|c|c|}
  \hline
  degree & singularities & (non-)split & type of result & reference\\
  \hline\hline
  $\ge 6$ & -- & split & asymptotic & \cite{MR1620682}\\
  \hline
  5 & -- & split & asymptotic & \cite{MR1909606}\\
  5 & -- & non-split & asymptotic & \cite{MR2099200}\\
  \hline
  4 & $\Dfive$ & split & asymptotic & 
  \cite{MR1906155}, \cite{math.NT/0412086}\\
  4 & $\Dfour$ & non-split & asymptotic & \cite{math.NT/0502510}\\
  4 & $\Dfour$ & split & asymptotic & this paper\\
  4 & $3\Aone$ & split & bounds & \cite{math.NT/0511041}\\
  \hline
  3 & $3\Atwo$ & split & asymptotic & 
  \cite{MR1620682}, \cite{MR2000b:11074}, \dots\\
  3 & $4\Aone$ & split & bounds & \cite{MR2075628}\\
  3 & $\Dfour$ & split & bounds & \cite{math.NT/0404245}\\
  3 & $\Esix$ & split & asymptotic & 
  \cite{math.NT/0504016}, \cite{math.NT/0509370}\\
  \hline
\end{tabular}
\smallskip
\caption{Results for Del Pezzo surfaces}
\label{tab:overview}
\end{table}

\begin{figure}[ht]
  \centering
  \includegraphics[width=11cm]{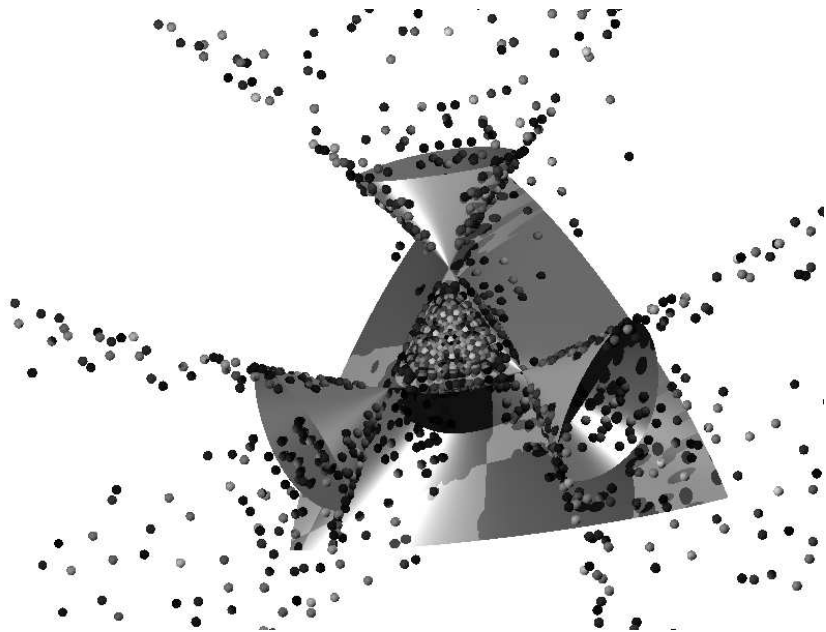}
  \caption{Points of height $\le 50$ on the Cayley cubic surface
    $x_0x_1x_2+x_0x_1x_3+x_0x_2x_3+x_1x_2x_3=0$.}
  \label{fig:cayley}
\end{figure}

\begin{figure}[ht]
  \centering
  \includegraphics[width=11cm]{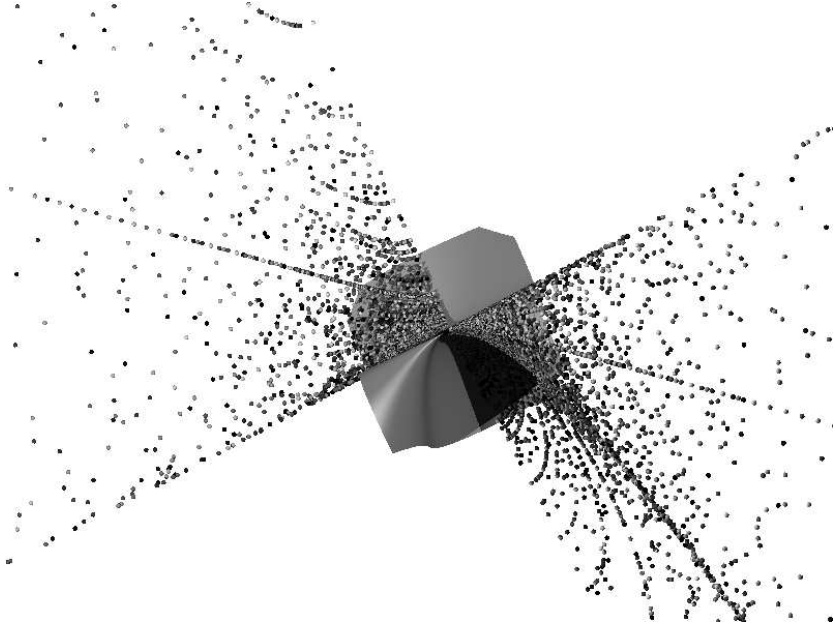}
  \caption{Points of height $\le 1000$ on the $\Esix$ singular cubic 
    surface $x_1x_2^2+x_2x_0^2+x_3^3=0$ with $x_0,x_2 > 0$.}
  \label{fig:e6}
\end{figure}

Figure \ref{fig:cayley} shows all points of height $\le 50$ on the
Cayley cubic surface (Example~\ref{ex:cayley}), which has four
singularities of type $\Aone$ and was considered in \cite{MR2075628}.
In Figure \ref{fig:e6}, we see points of height $\le 1000$ on the
$\Esix$ cubic surface (\cite{math.NT/0504016} and
\cite{math.NT/0509370}).

The proofs of Manin's conjecture proceed either via the height zeta
function
\[Z(s):=\sum_{\xx\in X^{\circ}(\QQ)} H_{-K_X}(\xx)^{-s},\]
whose analytic properties are related to the asymptotic
\eqref{eqn:asym} by Tauberian theorems, or via the lifting of the
counting problem to the \emph{universal torsor} -- an auxiliary
variety parametrizing rational points. The torsor approach has been
developed by Colliot-Th\'el\`ene and Sansuc in the context of the
Brauer-Manin obstruction \cite{MR89f:11082} and applied to Manin's
conjecture by Peyre \cite{MR1679842} and Salberger \cite{MR1679841}.

In the simplest case of hypersurfaces $X=X_f\subset \PP^n$ over $\QQ$,
with $n\ge 4$, this is exactly the passage from rational vectors
$\xx=(x_0,\dots, x_n)$, modulo the diagonal action of $\QQ^*$, to
primitive lattice points $(\ZZ_\prim^{n+1}\setminus 0)/\pm$. Geometrically,
we have
\[\begin{CD}
  \AA^{n+1} \setminus 0 @>\Gm>> \PP^n
\end{CD}
\qquad \text{and} \qquad 
\begin{CD}
  \TT_X @>\Gm>> X.
\end{CD}\] Here, $\TT_X$ is the hypersurface in $\AA^{n+1} \setminus
0$ defined by the form $f$, the 1-dimensional torus $\Gm$ is
interpreted as the N\'eron-Severi torus $T_{\NS}$, i.e., an algebraic
torus whose characters $\mathfrak X^*(T_{\NS})$ are isomorphic to the
Picard group (lattice) of $\PP^n$, resp.  $X$, and the map is the
natural quotient by its (diagonal) action.  Rational points on the
base are lifted to integral points on the torsor, modulo the action of
the group of units $T_{\NS}(\ZZ) =\{\pm 1\}$. The height inequality on
the base $H(\xx)\le B$ translates into the usual height inequality on
the torsor \eqref{eqn:height}.

In general, a torsor under an algebraic torus $T$ is determined by a
homomorpism $\chi\,:\, \mathfrak X^*(T)\rightarrow \Pic(X)$ to the
Picard group of the underlying variety $X$; the term \emph{universal}
is applied when $\chi$ is an isomorphism.

However, for hypersurfaces in $\Pthree$, or more generally for
complete intersection surfaces (i.e., $S$ is the intersection of $k$
hypersurfaces in $\PP^{k+2}$), the Picard group may have higher rank.
For example, for split smooth cubic surfaces $S=S_3\subset \Pthree$
the rank is 7, so that the dimension of the corresponding universal
torsor $\TT_S$ is 9; for quartic Del Pezzo surfaces these are 6 and 8,
respectively.

It is expected that the passage to universal torsors, which can be
considered as natural \emph{descent varieties}, will facilitate the
proof of Manin's conjecture (Conjecture~\ref{conj:manin}), at least
for Del Pezzo surfaces.  Rational points on $S$ are lifted to certain
integral points on $\TT_S$, modulo the action of $T_{\NS}(\ZZ)= (\pm
1)^{r}$, where $r$ is the rank of $\Pic(S)$, and the height inequality
on $S$ translates into appropriate inequalities on $\TT_S$.  This
explains the interest in the projective geometry of torsors, and
expecially, in their equations.  The explicit determination of these
equations is an interesting algebro-geometric problem, involving tools
from invariant theory and toric geometry.

In this note, we illustrate the torsor approach to asymptotics of
rational points in the case of a particular singular surface $S\subset \Pfour$ of
degree~4 given by:
\begin{equation}\label{eq:surface} 
x_0x_3-x_1x_4 = x_0x_1+x_1x_3+x_2^2 = 0.
\end{equation}
This is a split Del Pezzo surface, with a singularity of
type~$\Dfour$.

\begin{theorem}\label{thm:main}
  The number of $\QQ$-rational points of anticanonical height bounded
  by $B$ on the complement $\So$ of the $\QQ$-rational lines on $S$ as
  in \eqref{eq:surface} satisfies
  \[\NU = c_{S,H}\cdot B\cdot Q(\log B)+O(B (\log B)^3)
  \qquad\text{as $B \to \infty$},\]
  where $Q$ is a monic polynomial of degree $5$, and 
  \[c_{S,H}= \frac{1}{34560} \cdot \omega_\infty \cdot \prod_p
(1-1/p)^6(1+6/p+1/p^2)\] is the constant predicted by Peyre
\cite{MR1340296}, with $p$ running through all primes and
 \begin{equation*}
    \omega_\infty = 3\int\int\int_{\{(t,u,v)\in \RR^3 \mid 0\le v \le 1,
      |tv^2|,|v^2u|,|v(tv+u^2)|,|t(tv+u^2)|\le 1\}}1 \dd t \dd u \dd v.\
  \end{equation*}
\end{theorem}

In \cite{math.NT/0502510}, Manin's conjecture is
proved for a non-split surface with a singularity of the same type.
However, these results do not follow from each other.

In Section \ref{sec:geometry}, we collect some facts about the
geometric structure of $S$. In Section~\ref{sec:conformity}, we
calculate the expected value of $c_{S,H}$ and show that Theorem
\ref{thm:main} agrees with Manin's conjecture.

In our case, the universal torsor is an affine hypersurface.  In
Section~\ref{sec:torsor}, we calculate its equation, stressing the
relation with the geometry of $S$.  We make explicit the coprimality
and the height conditions.  The method is more systematic than the
derivation of torsor equations in \cite{math.NT/0412086} and
\cite{math.NT/0509370}, and should bootstrap to more complicated
cases, e.g., other split Del Pezzo surfaces.

Note that our method gives coprimality conditions which are different
from the ones in \cite{math.NT/0412086} and \cite{math.NT/0509370},
but which are in a certain sense more natural: they are related to the
set of points on $\TT_S$ which are \emph{stable} with respect to the
action of the N\'eron-Severi torus (in the sense of geometric
invariant theory, c.f., \cite{MR2004511} and
\cite{MR2001i:14059}). Our conditions involve only coprimality of
certain pairs of variables, while the other method produces a mix of
square-free variables and coprimalities.

In Section \ref{sec:summations}, we estimate the number of integral
points on the universal torsor by iterating summations over the torsor
variables and using results of elementary analytic number theory.
Finally we arrive at Lemma \ref{lem:sum_eta}, which is very similar in
appearance to \cite[Lemma 10]{math.NT/0412086} and \cite[Lemma
12]{math.NT/0504016}.  In Section \ref{sec:proof} we use familiar
methods of height zeta functions to derive the exact asymptotic.  We
isolate the expected constant $c_{S,H}$ and finish the proof of
Theorem \ref{thm:main}.  In Section~\ref{sec:torsorexamples} we write down
examples of universal torsors for other Del Pezzo surfaces and discuss
their geometry.

\

\noindent{\bf Acknowledgment.}
Part of this work was done while the authors were visiting the CRM at
the Universit\'e de Montr\'eal during the special year on \emph{Analysis in
Number Theory}. We are grateful for the invitation and ideal working
conditions.

\section{Geometric background}\label{sec:geometry}

In this section, we collect some geometric facts concerning the
surface $S$. We show that Manin's conjecture for $S$ is not a special
case of available more general results for Del Pezzo surfaces.

\begin{lemma}\label{lem:geometry}
  The surface $S$ has the following properties:
  \begin{enumerate}
  \item \label{item:geometry_singularity} It has exactly one
    singularity of type $\Dfour$ at $q = (0:0:0:0:1)$.
  \item \label{item:lines} $S$ contains exactly two lines: \[E_5 =
    \{x_0=x_1=x_2=0\}\,\,\text{ and }\,\, E_6 = \{x_1=x_2=x_3=0\},\]
    which intersect in $q$.
  \item \label{item:geometry_projection} The projection from the line
    $E_5$ is a birational map 
    \[\begin{array}{cccc}
      \phi: & S & \rto & \Ptwo\\
      & \xx & \mapsto & (x_0:x_2:x_1)
    \end{array}\]
    which is defined outside $E_5$. It restricts to an
    isomorphism between 
    \[\So = S \setminus (E_5 \cup E_6) = \{\xx \in S
    \mid x_1 \ne 0\}\,\text{ and }\, \AA^2 \cong \{(t:u:v) \mid v
    \ne 0\} \subset \Ptwo,\] whose inverse is the restriction of
    \[\begin{array}{cccl}
      \psi : & \Ptwo & \rto & S,\\
      & (t:u:v) & \mapsto & (tv^2:v^3:v^2u:-v(tv+u^2):-t(tv+u^2))
    \end{array}\]
    Similar results hold for the projection from $E_6$.
  \item \label{item:geometry_dynkin} The process of resolving the
    singularity $q$ gives four exceptional curves $E_1, \dots, E_4$
    and produces the minimal desingularization $\tS$, which is also
    the blow-up of $\Ptwo$ in five points.
  \end{enumerate}
\end{lemma}

\begin{proof}
Direct computations.
\end{proof}

It will be important to know the details of the sequence of five
blow-ups of $\Ptwo$ giving $\tS$ as in Lemma
\ref{lem:geometry}(\ref{item:geometry_dynkin}):

In order to describe the points in $\Ptwo$, we need the lines 
\[E_3 = \{v=0\},\qquad A_1 = \{u=0\},\qquad A_2 = \{t=0\}\]
and the curve $A_3 = \{tv+u^2=0\}$.

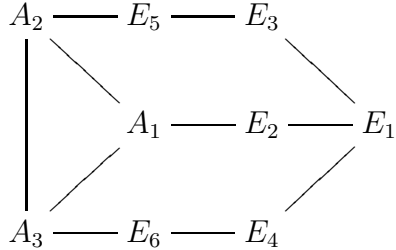
\begin{figure}[ht]
  \centering
  \[\xymatrix{A_2  \ar@{-}[r]\ar@{-}[dr] \ar@{-}[dd]& E_5 \ar@{-}[r] 
    & E_3 \ar@{-}[dr]\\
    & A_1 \ar@{-}[r] & E_2 \ar@{-}[r] & E_1 \\
    A_3 \ar@{-}[r] \ar@{-}[ur] & E_6 \ar@{-}[r] & E_4 \ar@{-}[ur]}\]
  \caption{Extended Dynkin diagram}
  \label{fig:dynkin}
\end{figure}

\begin{lemma}\label{lem:blow-ups}
  The following five blow-ups of $\Ptwo$ result in $\tS$:
  \begin{itemize}
  \item Blow up the intersection of $E_3, A_1, A_3$, giving $E_2$.
  \item Blow up the intersection of $E_2, E_3, A_3$, giving $E_1$.
  \item Blow up the intersection of $E_1$ and $A_3$, giving $E_4$.
  \item Blow up the intersection of $E_4$ and $A_3$, giving $E_6$.
  \item Blow up the intersection of $E_3$ and $A_2$, giving $E_5$.
  \end{itemize}
  Here, the order of the first four blow-ups is fixed, and the fifth
  blow-up can be done at any time.
  
  The Dynkin diagram in Figure \ref{fig:dynkin} describes the final
  configuration of divisors $E_1, \dots, E_6, A_1, A_2, A_3$.  Here,
  $A_1, A_2, A_3$ intersect at one point.
\end{lemma}

The quartic Del Pezzo surface with a singularity of type $\Dfour$ is
not toric, and Manin's conjecture does not follow from the results of
\cite{MR1620682}.  The $\Dfive$ example of \cite{math.NT/0412086} is
an equivariant compactification of $\Ga^2$ (i.e., $S$ has a Zariski open
subset isomorphic to $\AA^2$, and the obvious action of $\Ga^2$ on
this open subset extends to $S$), and thus a special case of
\cite{MR1906155}.

\begin{lemma}\label{lem:compact}
  The quartic Del Pezzo surface with a singularity of type $\Dfour$ is a
  compactification of $\AA^2$, but not an equivariant compactification of
  $\Ga^2$.
\end{lemma}

\begin{proof}
  We follow the strategy of \cite[Remark 3.3]{MR2029868}.
  
  Consider the maps $\phi, \psi$ as in Lemma
  \ref{lem:geometry}(\ref{item:geometry_projection}). As $\psi$
  restricts to an isomorphism between $\AA^2$ and the open set $\So
  \subset S$, the surface $S$ is a compactification of $\AA^2$.
  
  If $S$ were an equivariant compactification of $\Ga^2$ then the
  projection $\phi$ from $E_5$ would be a $\Ga^2$-equivariant map,
  giving a $\Ga^2$-action on $\Ptwo$.  The line $\{v = 0\}$ would be
  invariant under this action. The only such action is the standard
  translation action
  \[\begin{array}{cccc}
    \tau:&\Ptwo &\to &\Ptwo,\\
    &(t:u:v) &\mapsto &(t+\alpha v: u+\beta v: v).
  \end{array}\]
  However, this action does not leave the linear series
  \[(tv^2:v^3:v^2u:-v(tv+u^2):-t(tv+u^2))\] invariant, 
  which can be seen after calculating
  \begin{equation*}
    \begin{split}
      t(tv+u^2) \mapsto &(t+\alpha v)((t+\alpha v)v+(u+\beta v)^2)\\
      = &t(tv+u^2) + 2\beta tuv + (\beta^2+\alpha) tv^2 + \alpha
      v(tv+u^2)\\ &+ 2\alpha\beta v^2u +
      (\alpha\beta^2+\alpha^2) v^3,
    \end{split}
  \end{equation*}
  since the term $tuv$ does not appear in the original linear series.
\end{proof}

\section{Manin's conjecture}\label{sec:conformity}

\begin{lemma}\label{lem:conjecture}
Let $S$  be the surface \eqref{eq:surface}. Manin's conjecture
for $S$ states that the number of rational points of 
height $\le B$ outside the two lines is given by 
\[\NU \sim c_{S,H}\cdot B(\log B)^5,\]
  where $c_{S,H} = \alpha(S)\cdot \beta(S)\cdot \omega_H(S)$ with 
  \begin{equation*}
    \begin{split}
      \alpha(S) &= (5!\cdot 4\cdot 2\cdot 3 \cdot 3 \cdot 2 \cdot
      2)^{-1} = (34560)^{-1}\\
      \beta(S) &= 1\\
      \omega_H(S) &= \omega_\infty \cdot \prod_p (1-1/p)^6(1+6/p+1/p^2)
    \end{split}
  \end{equation*}
  and
  \begin{equation*}
    \omega_\infty = 3\int\int\int_{\{(t,u,v)\in \RR^3 \mid 0\le v \le 1,
      |tv^2|,|v^2u|,|v(tv+u^2)|,|t(tv+u^2)|\le 1\}}1 \dd t \dd u \dd v\
  \end{equation*}
\end{lemma}

\begin{proof}
  Since $S$ is split over $\QQ$, we have $\rk(\Pic(\tS)) = 6$, and the
  expected exponent of $\log B$ is $5$.  Further, $\beta(S) = 1$. The
  computation of $c_{S,H}$ is done on the desingularization $\tS$. For
  the computation of $\alpha(S)$, observe that the effective cone of
  $\tS$ in $\Pic(\tS)$ is simplicial (it is generated by the
  exceptional curves $E_1, \dots, E_6$, and their number equals the
  rank of $\Pic(\tS)$), and
  \[\anti = 4E_1+2E_2+3E_3+3E_4+2E_5+2E_6.\]
  The calculation is analog to \cite[Lemma 2]{math.NT/0504016} (see
  \cite{del_pezzo_alpha} for its calculation in general).  The
  constant $\omega_H(S)$ is computed as in
  \cite[Lemma~1]{math.NT/0412086} and \cite[Lemma~2]{math.NT/0504016}.
\end{proof}

\section{The universal torsor}\label{sec:torsor}

As explained above, the problem of counting rational points of bounded
height on the surface $S$ translates into a counting problem for
certain integral points on the universal torsor, subject to
coprimality and height inequalities.  In the first part of this
section, we describe these conditions in detail.  They are obtained by
a process of introducing new variables which are the greatest common
divisors of other variables.  Geometrically, this corresponds to the
realization of $\tS$ as a blow-up of $\Ptwo$ in five points.

In the second part, we prove our claims.

The universal torsor $\TT_S$ of $S$ is an open subset of the
hypersurface in $\AA^9 = \Spec \ZZ[\eta_1, \dots,
\eta_6,\alpha_1,\alpha_2,\alpha_3]$ defined by the equation
\begin{equation}\label{eq:torsor}
T(\ee,\aa) = \Tone + \Ttwo + \Tthr = 0.
\end{equation}
The projection $\Psi:\TT_S\to S$ is defined by 
\begin{equation}\label{eq:projection}
  (\Psi^*(x_i)) = (\sbase 2 1 2 1 2 0 \atwo, \sbase 4 2 3 3 2 2,
  \sbase 3 2 2 2 1 1 \aone, \sbase 2 1 1 2 0 2 \athr, \atwo\athr),
\end{equation}
where we use the notation $\sbase {n_1}{n_2}{n_3}{n_4}{n_5}{n_6} =
\eta_1^{n_1}\eta_2^{n_2}\eta_3^{n_3}\eta_4^{n_4}\eta_5^{n_5}\eta_6^{n_6}$.

The coprimality conditions can be derived from the extended Dynkin
diagram (see Figure~\ref{fig:dynkin}). Two variables are allowed to
have a common factor if and only if the corresponding divisors ($E_i$
for $\eta_i$ and $A_i$ for $\ai$) intersect (i.e., are connected by an
edge in the diagram). Furthermore, $\gcd(\alpha_1,\alpha_2,\alpha_3) >
1$ is allowed (corresponding to the fact that $A_1, A_2, A_3$
intersect in one point).

We will show below that there is a bijection between rational points
on $\So \subset S$ and integral points on an open subset of $\TT_S$,
subject to these coprimality conditions.

We will later refer to
\begin{align}
  \label{eq:coprim_eta}&\text{coprimality between $\eta_i$ as in Figure~\ref{fig:dynkin},}\\
  \label{eq:coprim_a1}&\gcd(\aone,\Cone)=1,\\
  \label{eq:coprim_a2}&\gcd(\atwo,\Ctwo)=1,\\
  \label{eq:coprim_a3}&\gcd(\athr,\Cthr)=1.
\end{align}

To count the number of $\xx \in S(\QQ)$ such that $H(\xx) \le B$, we
must lift this condition to the universal torsor, i.e.,
$H(\Psi(\ee,\aa)) \le B$. This is the same as \[|\sbase 2 1 2 1 2 0
\atwo| \le B, \qquad\dots,\qquad |\atwo\athr|\le B,\] using the five
monomials occuring in \eqref{eq:projection}. These have no common
factors, provided the coprimality conditions are fulfilled (direct
verification).

It will be useful to write the height conditions as follows: let
\[X_0 = \left(\frac{\sbase 4 2 3 3 2 2}{B}\right)^{1/3},\, X_1 = (B\sbase {-1}{-2}0 0 1
1)^{1/3},\, X_2 = (B\sbase 2 1 0 3 {-2} 4)^{1/3}.\]
Then 
\begin{align}
  \label{eq:height_eta}&|X_0^3| \le 1\\
  \label{eq:height_a1}&|X_0^2 (\alpha_1/X_1)| \le 1\\
  \label{eq:height_a2}
  \begin{split}
    &|X_0^2 (\alpha_2/X_2)| \le 1, \qquad
    |X_0(X_0(\alpha_2/X_2)+(\alpha_1/X_1)^2)|\le 1, \qquad\\
    &|(\alpha_2/X_2)(X_0(\alpha_2/X_2)+(\alpha_1/X_1)^2)| \le 1
\end{split}
\end{align}
are equivalent to the five height conditions. Here we have used the
torsor equation to eliminate $\athr$ because in our counting argument
we will also use that $\athr$ is determined by the other variables.

We now prove the above claims.

\begin{lemma}\label{lem:bijection}
  The map $\Psi$ gives a bijection between the set of points $\xx$ of
  $\So(\QQ)$ such that $H(\xx) \le B$ and the set 
  \[\TT_1 := \Bigg\{(\ee,\aa)\in \ZZp^6\times\ZZ^3 \Bigg|
  \begin{aligned}&\text{equation} ~\eqref{eq:torsor}, \\ 
                 &\text{coprimality} ~ 
      \eqref{eq:coprim_eta}, \eqref{eq:coprim_a1},
      \eqref{eq:coprim_a2}, \eqref{eq:coprim_a3},\\
    &\text{inequalities \eqref{eq:height_eta}, \eqref{eq:height_a1},
      \eqref{eq:coprim_a3} hold}
\end{aligned}
\Bigg\}\]
\end{lemma}

\begin{proof}
  The map $\psi$ of
  Lemma~\ref{lem:geometry}(\ref{item:geometry_projection}) induces a
  bijection \[\psi_0: (\eta_3, \aone, \atwo) \mapsto (\eta_3^2\atwo,
  \eta_3^3, \eta_3^2\aone, \eta_3\athr, \atwo\athr),\] where $\athr
  :=-(\eta_3\atwo+\aone^2)$, i.e., \[T_0 := \aone^2+\eta_3\atwo+\athr =
  0,\] between
  \[\{(\eta_3,\aone, \atwo) \in \ZZp \times \ZZ^2 \mid \gcd(\eta_3,
  \aone, \atwo) = 1\} 
\,\, \text{ and }\,\, \So(\QQ) \subset S(\QQ).
\]
  The height function on $\So(\QQ)$ is given by\[H(\psi_0(\eta_3, \aone,
  \atwo)) = \frac{\max(|\eta_3^2\atwo|, |\eta_3^3|, |\eta_3^2\aone|,
    |\eta_3\athr|, |\atwo\athr|)}{\gcd(\eta_3^2\atwo, \eta_3^3,
    \eta_3^2\aone, \eta_3\athr, \atwo\athr)}.\]
  
  The derivation of the torsor equation from the map $\psi_0$ together
  with the coprimality conditions and the lifted height function is
  parallel to the blow-up process described in
  Lemma~\ref{lem:blow-ups}. More precisely, each line $E_3, A_1, A_2$
  in $\Ptwo$ corresponds to a coordinate function $\eta_3, \aone,
  \atwo$ vanishing in one of the lines; the blow-up of the
  intersection of two divisors gives an exceptional curve $E_i$,
  corresponding to the introduction of a new variable $\eta_i$ as the
  greatest common divisor of two old variables. Two divisors are
  disjoint if and only if the corresponding variables are coprime.
  This is summarized in Table~\ref{tab:dictionary}.
  \begin{table}[htbp]\centering
    \begin{tabular}{|l||l|}
      \hline
      Variables, Equations & Geometry \\
      \hline\hline
      variables & divisors\\
      initial variables & coordinate lines\\
      $\eta_3,\aone, \atwo$ & $E_3,A_1, A_2$\\
      \hline
      taking $\gcd$ of two variables & blowing up intersection of divisors\\
      new $\gcd$-variable & exceptional curve\\
      $\eta_2, \eta_1, \eta_4, \eta_6, \eta_5$ & $E_2, E_1, E_4, E_6, E_5$\\
      \hline
      extra variable & extra curve\\
      $\athr$ & $A_3$\\
      starting relation & starting description\\
      $\athr=-(\eta_3\atwo+\aone^2)$ & $A_3 = \{\eta_3\atwo+\aone^2=0\}$\\
      final relation & torsor equation\\
      $\Tthr = -(\Ttwo+\Tone)$ & $\Tone+\Ttwo+\Tthr=0$\\
      \hline
    \end{tabular}
    \smallskip
    \caption{Dictionary between $\gcd$-process and blow-ups}
    \label{tab:dictionary}
  \end{table}
  
  This plan will now be implemented in five steps; at each step, the
  map \[\psi_i: \ZZp^{i+1} \times \ZZ^3 \to \So(\QQ)\] gives a bijection
  between:
  \begin{itemize}
  \item the set of all $(\eta_j, \aone, \atwo, \athr) \in \ZZp^{i+}
    \times \ZZ^3$ satisfying certain coprimality conditions (described
    by the extended Dynkin diagram corresponding to the $i$-th blow-up
    of Lemma~\ref{lem:blow-ups}), an equation $T_i$,
    \[H(\psi_i(\eta_j, \alpha_j)) = \frac{\max_k(|\psi_i(\eta_j,
      \alpha_j)_k|)}{\gcd(\psi_i(\eta_j, \alpha_j)_k)} \le B,\]
  \item the set of all $\xx \in \So(\QQ)$ with $H(\xx) \le B$.
  \end{itemize}
  The steps are as follows:
  \begin{enumerate}
  \item Let $\eta_2 := \gcd(\eta_3, \aone) \in \ZZp$. Then \[\eta_3 =
    \eta_2\eta_3',\qquad \aone = \eta_2\aone',\qquad\text{ with
      $\gcd(\eta_3', \aone') = 1$.}\] Since $\eta_2 \mid \athr$, we
    can write $\athr = \eta_2\athr'$.  Then $\athr' =
    -(\eta_3'\atwo+\eta_2\aone'^2)$. After renaming the variables, we
    have \[T_1 = \eta_2\aone^2 + \eta_3\atwo + \athr = 0\] and
    \[\psi_1: (\eta_2, \eta_3, \aone, \atwo, \athr) \mapsto 
    (\eta_2\eta_3^2\atwo: \eta_2^2\eta_3^3: \eta_2^2\eta_3^2\aone:
    \eta_2\eta_3\athr: \atwo\athr).\] Here, we have eliminated the
    common factor $\eta_2$ which occured in all five components of the
    image. Below, we repeat the corresponding transformation at each
    step.
  \item Let $\eta_1 := \gcd(\eta_2, \eta_3) \in \ZZp$. Then \[\eta_2 =
    \eta_1\eta_2', \qquad \eta_3 = \eta_1\eta_3',\qquad \text{with
      $\gcd(\eta_2', \eta_3')=1$.}\] As $\eta_1 \mid \athr$, we write
    $\athr = \eta_1\athr'$, and we obtain: \[T_2 = \eta_2\aone^2 +
    \eta_3\atwo + \athr = 0\] and
    \begin{multline*}
    \psi_2: (\eta_1, \eta_2, \eta_3, \aone, \atwo, \athr) 
    \mapsto \\
    (\eta_1^2\eta_2\eta_3^2\atwo : \eta_1^4\eta_2^2\eta_3^3 :
    \eta_1^3\eta_2^2\eta_3^2\aone : \eta_1^2\eta_2\eta_3\athr :
    \atwo\athr). 
    \end{multline*}
  \item Let $\eta_4 := \gcd(\eta_1, \athr) \in \ZZp$. Then \[\eta_1 =
    \eta_4\eta_1',\qquad \athr = \eta_4\athr',\qquad \text{with
      $\gcd(\eta_1', \athr')=1$.}\] We get after removing $'$ again:
    \[T_3 = \eta_2\aone^2 + \eta_3\atwo + \eta_4\athr = 0\] and
    \begin{multline*}
      \psi_3 : (\eta_1, \eta_2, \eta_3, \eta_4, \aone, \atwo, \athr)
      \mapsto\\
      (\eta_1^2\eta_2\eta_3^2\eta_4\atwo :
      \eta_1^4\eta_2^2\eta_3^3\eta_4^3 :
      \eta_1^3\eta_2^2\eta_3^2\eta_4^2\aone :
      \eta_1^2\eta_2\eta_3\eta_4^2\athr : \atwo\athr).
    \end{multline*}
  \item Let $\eta_6 := \gcd(\eta_4, \athr) \in \ZZp$. Then \[\eta_4 =
    \eta_6\eta_4',\qquad \athr = \eta_6\athr', \qquad\text{with
      $\gcd(\eta_4', \athr') = 1$.}\] We obtain \[T_4 = \eta_2\aone^2 +
    \eta_3\atwo + \eta_4\eta_6^2\athr = 0\] and
    \begin{multline*}
      \psi_4 :(\eta_1, \eta_2, \eta_3, \eta_4, \eta_6, 
      \aone, \atwo, \athr) \mapsto \\
      (\eta_1^2\eta_2\eta_3^2\eta_4\atwo :
      \eta_1^4\eta_2^2\eta_3^3\eta_4^3\eta_6^2 :
      \eta_1^3\eta_2^2\eta_3^2\eta_4^2\eta_6\aone :
      \eta_1^2\eta_2\eta_3\eta_4^2\eta_6^2\athr : \atwo\athr).
    \end{multline*}
  \item The final step is $\eta_5 := \gcd(\eta_3, \atwo) \in \ZZp$, we 
    could have done it earlier (just as the blow-up of
    the intersection of $E_3, A_2$ in Lemma \eqref{eq:projection}). Then
    \[\eta_3 = \eta_5\eta_3',\qquad \atwo = \eta_5\atwo',\qquad\text{with
      $\gcd(\eta_3', \atwo')=1.$}\]  We get \[T_5 = \eta_2\aone^2 +
    \eta_3\eta_5\atwo + \eta_4\eta_6^2\athr = 0\] and
    \begin{multline*}
      \psi_5 :(\eta_1, \eta_2, \eta_3, \eta_4, \eta_5, \eta_6, \aone,
      \atwo, \athr) \mapsto\\
      (\eta_1^2\eta_2\eta_3^2\eta_4\eta_5^2\atwo :
      \eta_1^4\eta_2^2\eta_3^3\eta_4^3\eta_5^2\eta_6^2 :
      \eta_1^3\eta_2^2\eta_3^2\eta_4^2\eta_5\eta_6\aone :
      \eta_1^2\eta_2\eta_3\eta_4^2\eta_6^2\athr : \atwo\athr)
    \end{multline*}
  \end{enumerate}
  We observe that at each stage 
  the coprimality conditions
  correspond to intersection properties of the respective divisors. 
  The final result is summarized in
  Figure~\ref{fig:dynkin}, which encodes data from \eqref{eq:coprim_eta},
  \eqref{eq:coprim_a1}, \eqref{eq:coprim_a2}, \eqref{eq:coprim_a3}.
  
  Note that $\psi_5$ is  $\Psi$ from \eqref{eq:projection}.
  As mentioned above, $\gcd(\psi_5(\eta_j, \alpha_j)_k)$ 
  (over all five components of the image) is trivial by the coprimality
  conditions of Figure~\ref{fig:dynkin}. Therefore, 
  $H(\psi_5(\ee, \aa)) \le B$ is equivalent to \eqref{eq:height_eta},
  \eqref{eq:height_a1}, \eqref{eq:height_a2}.
  
  Finally, $T_5$ is the torsor equation $T$ \eqref{eq:torsor}.
\end{proof}

\section{Summations}\label{sec:summations}

In the first step, we estimate the number of $(\alpha_1, \alpha_2,
\alpha_3) \in \ZZ^3$ which fulfill the torsor equation $T$
\eqref{eq:torsor} and the height and coprimality conditions. For fixed
$(\alpha_1, \alpha_2)$, the torsor equation $T$ has a
solution $\alpha_3$ if and only if the congruence
\[\congr{\Tone+\Ttwo}0\Fthr\] holds and the conditions on the height and
coprimalities are fulfilled.

We have already written the height conditions so that they do not
depend on $\athr$. For the coprimality, we must ensure that
\eqref{eq:coprim_a2} and \eqref{eq:coprim_a3} are fulfilled.

As $\gcd(\Ftwo, \Fthr)=1$, we can find the multiplicative inverse
$c_1$ of $\Ftwo$ modulo $\Fthr$, so that
\begin{equation}
\label{eq:c1c2}
c_1\Ftwo = 1+c_2\Fthr
\end{equation}
for a suitable $c_2$.
Choosing 
\begin{align}
  \label{eq:c1c3a2}\atwo &= c_3\Fthr - c_1\Tone,\\
  \label{eq:c2c3a3}\athr &= c_2\Tone-c_3\Ftwo
\end{align}
gives a solution of \eqref{eq:torsor} for any $c_3 \in \ZZ$.

Without the coprimality conditions, the number of pairs
$(\alpha_2,\alpha_3)$ satisfying $T$ and \eqref{eq:height_a2} would
differ at most by $O(1)$ from $1/\Fthr$ of the length of the interval
described by \eqref{eq:height_a2}.  However, the coprimality
conditions \eqref{eq:coprim_a2} and \eqref{eq:coprim_a3} impose
further restrictions on the choice of $c_3$.  A slight complication
arises from the fact that because of $T$, some of the conditions are
fulfilled automatically once $\ee, \aone$ satisfy
\eqref{eq:coprim_eta} and \eqref{eq:coprim_a1}.

Conditions~\eqref{eq:coprim_eta} imply that the possibilities for a
prime $p$ to divide more than one of the $\eta_i$ are very limited. We
distinguish twelve cases, listed in Column~2 of
Table~\ref{tab:coprim_alpha}.

\begin{table}[ht]
  \centering
  \[\begin{array}{|c|c||c|c|c|}
    \hline
    \text{case} & p\mid \dots & p\mid \alpha_1 & p\mid \alpha_2 & 
    p\mid \alpha_3\\
    \hline\hline
    0 & - & \allow & \allow & \allow\\
    i & \eta_1 & \restr & \restr & \restr\\
    ii & \eta_2 & \allow & \restr & \autom\\
    iii & \eta_3 & \restr & \restr & \autom\\
    iv & \eta_4 & \restr & \autom & \restr\\
    v & \eta_5 & \restr & \allow & \autom\\
    vi & \eta_6 & \restr & \autom & \allow\\
    \hline
    vii & \eta_1,\eta_2 & \restr & \restr & \autom\\
    viii & \eta_1,\eta_3 & \restr & \restr & \autom\\
    ix & \eta_1,\eta_4 & \restr & \autom & \restr\\
    x & \eta_3,\eta_5 & \restr & \restr & \autom\\
    xi & \eta_4,\eta_6 & \restr & \autom & \restr\\
    \hline
  \end{array}\]
  \smallskip
  \caption{Coprimality conditions}
  \label{tab:coprim_alpha}
\end{table}

In Columns 4 and 5, we have denoted the relevant information for the
divisibility of $\atwo, \athr$ by primes $p$ which are divisors of the
$\eta_i$ in Column 2, but of no other $\eta_j$:
\begin{itemize}
\item ``allowed'' means that $\ai$ may be divisible by $p$.
\item ``automatically'' means that the conditions on the $\eta_i$ and
  the other $\alpha_j$ imply that $p \nd \ai$.  These
  two cases do not impose conditions on $c_3$ modulo $p$.
\item ``restriction'' means that $c_3$ is not allowed to be in a
  certain congruence class modulo $p$ in order to fulfill the
  condition that $p$ must not divide $\ai$.
\end{itemize}

The information in the table is derived as follows:
\begin{itemize}
\item If $p\mid \eta_3$, then $p \nd c_2$ from \eqref{eq:c1c2}, and $p
  \nd \alpha_1\eta_2$ because of \eqref{eq:coprim_eta},
  \eqref{eq:coprim_a1}, so by \eqref{eq:c2c3a3}, $p \nd \athr$
  independently of the choice of $c_3$. Since $p \nd \Fthr$, we see
  from \eqref{eq:c1c3a2} that $p \mid \atwo$ for one in $p$ subsequent
  choices of $c_3$ which we must therefore exclude. This explains
  cases $iii$ and $viii$.
\item In case $vii$, the same is true for $\atwo$. More precisely, we
  see that we must exclude $\congr{c_3}0 p$. By \eqref{eq:c2c3a3}, $p
  \nd c_3$ implies that $p\nd \athr$, so we do not need another
  condition on $c_3$.
\item In case $i$, we see that $p\mid \atwo$ for one in $p$ subsequent
  choices of $c_3$, and the same holds for $\athr$. However, in this
  case, $p$ cannot divide $\atwo, \athr$ for the same choice of $c_3$,
  as we can see by considering $T$: since $p \nd \Tone$, it is
  impossible that $p\mid \atwo,\athr$. Therefore, we must exclude two
  out of $p$ subsequent choices of $p$ in order to fulfill $p \nd
  \atwo,\athr$.
\item In the other cases, the arguments are similar.
\end{itemize}

The number of $(\atwo, \athr) \in \ZZ^2$ subject to $T$,
\eqref{eq:coprim_a2}, \eqref{eq:coprim_a3}, \eqref{eq:height_a2}
equals the number of $c_3$ such that $\atwo, \athr$ as in
\eqref{eq:c1c3a2}, \eqref{eq:c2c3a3} satisfy these conditions. This
can be estimated as $1/\Fthr$ of the interval described by
\eqref{eq:height_a2}, multiplied by a product of local factors whose
value can be read off from Columns 2, 4, 5 of
Table~\ref{tab:coprim_alpha}: the divisibility properties of $\eta_i$
by $p$ determine whether zero, one or two out of $p$ subsequent values
of $c_3$ have to be excluded. Different primes can be considered
separately, and we define
\[\vartheta_{1,p} := 
\begin{cases}
  1-2/p, &\text{case $i$,} \\
  1-1/p, &\text{cases $ii-iv, vi-xi$,}\\
  1, &\text{case $0,v$.}
\end{cases}
\]
Let \[\vartheta_1(\ee) = \prod_p \vartheta_{1,p}\] be the product of
these local factors, and
\begin{equation}\label{eq:g_1}
g_1(u,v) = \int_{\{t \in \RR \mid |tv^2|,|t(tv+u^2)|,|v(tv+u^2)|\le 1\}} 
1 \dd t.
\end{equation}
Let $\omega(n)$ denote the number of primes dividing $n$.

\begin{lemma}\label{lem:sum_a2a3}
  For fixed $(\ee, \aone) \in \ZZp^6 \times \ZZ$ as in
  \eqref{eq:coprim_eta}, \eqref{eq:coprim_a1}, \eqref{eq:height_eta},
  \eqref{eq:height_a1}, the number of $(\atwo,\athr) \in \ZZ^2$
  satisfying $T$, \eqref{eq:coprim_a2}, \eqref{eq:coprim_a3},
  \eqref{eq:height_a2} is \[\NN_1(\ee,\aone) =
  \frac{\vartheta_1(\ee)X_2}{\Fthr} g_1(\aone/X_1,X_0) +
  O(2^{\omega(\Ctwo)}).\] The sum of error terms for all possible
  values of $(\ee, \aone)$ is $\ll B(\log B)^3$.
\end{lemma}

\begin{proof}
  The number of $c_3$ such that the resulting $\atwo, \athr$ satisfy
  \eqref{eq:height_a2} differs from
  $\frac{X_2}{\Fthr}g_1(\aone/X_1,X_0)$ by at most $O(1)$.
  
  Each $\vartheta_{1,p} \ne 1$ corresponds to a congruence condition
  on $c_3$ imposed by one of the cases $i-iv,vi-xi$. For each
  congruence condition, the actual ratio of allowed $c_3$ can differ
  at most by $O(1)$ from the $\vartheta_{1,p}$. The total number of
  these primes $p$ is \[\omega(\Ctwo) \ll 2^{\omega(\Ctwo)},\] which
  is independent of $\eta_5$ since any prime dividing only $\eta_5$
  contributes a trivial factor (see case $v$).
  
  Using the estimate \eqref{eq:height_a1} for $\aone$ in the first
  step and ignoring \eqref{eq:coprim_eta} \eqref{eq:coprim_a1}, which
  can only increase the error term, we obtain:
  \begin{equation*}
      \sum_{\ee}\sum_{\aone}
      2^{\omega(\Ctwo)} \le \sum_{\ee} \frac {B\cdot 2^{\omega(\Ctwo)}}
      {\sbase 3 2 2 2 1 1} \ll B(\log B)^3.
  \end{equation*}
  Here, we use $2^{\omega(n)} \lle n^\ep$ for the summations over
  $\eta_1, \eta_2, \eta_3, \eta_4$. For $\eta_6$, we employ
  \[\sum_{n\le x} 2^{\omega(n)} \ll x (\log x)\] together with partial
  summation, contributing a factor $(\log B)^2$, while the summation
  over $\eta_5$ gives another factor $\log B$.
\end{proof}

Next, we sum over all $\alpha_1$ subject to the coprimality condition
\eqref{eq:coprim_a1} and the height condition \eqref{eq:height_a1}.
Let
\begin{equation}\label{eq:g_2}
g_2(v) = \int_{\{u \in \RR \mid |v^2u| \le 1\}} g_1(u,v) \dd u
\end{equation}
Similar to our discussion for $\atwo, \athr$, the number of possible
values for $\aone$ as in \eqref{eq:height_a1}, while ignoring
\eqref{eq:coprim_a1} for the moment, is $X_1g_2(X_0) + O(1)$.

None of the coprimality conditions are fulfilled automatically, and
only common factors with $\eta_2$ are allowed (see Column 3 of
Table~\ref{tab:coprim_alpha}). Therefore, each prime factor of $\Cone$
reduces the number of allowed $\aone$ by a factor of $\vartheta_{2,p}
= 1-1/p$ with an error of at most $O(1)$. For all other primes $p$,
let $\vartheta_{2,p} = 1$, and let
\[\vartheta_2(\ee) = \prod_p \vartheta_{2,p}\qquad\text{and}\qquad 
\vartheta(\ee) =
\begin{cases}
  \vartheta_1(\ee)\cdot\vartheta_2(\ee), 
  &\text{\eqref{eq:coprim_eta} holds}\\
  0, &\text{otherwise.}
\end{cases}\]

\begin{lemma}\label{lem:sum_a1}
  For fixed $\eta \in \ZZp^6$ as in \eqref{eq:coprim_eta},
  \eqref{eq:height_eta}, the sum of $\NN_1(\ee, \aone)$ over all
  $\aone \in \ZZ$ satisfying \eqref{eq:coprim_a1}, \eqref{eq:height_a1}
  is
  \[\NN_2(\ee) := \frac{\vartheta(\ee)X_1X_2}{\Fthr} g_2(X_0) 
  + \R_2(\ee),\] where the sum of error terms $\R_2(\ee)$ over all
  possible $\ee$ is $\ll B\log B$.
\end{lemma}

\begin{proof}
  Let \[\NN(b_1,b_2) = \vartheta_1(\ee) \cdot \#\{\aone \in [b_1,b_2]
  \mid \gcd(\aone,\Cone)=1\}.\] Using M\"obius inversion, this is
  estimated as
  \[\NN(b_1,b_2) = \vartheta_1(\ee) \cdot \vartheta_2(\ee) \cdot (b_2-b_1) + 
  \R(b_1,b_2)\] with $\R(b_1,b_2) = O(2^{\omega(\Cone)})$. By partial
  summation, 
  \[\NN_2(\ee) =
  \frac{\vartheta(\ee)X_1X_2}{\Fthr} g_2(X_0) + \R_2(\ee)\] with
  \[\R_2(\ee) = \frac{-X_2}{\Fthr}\int_{\{u \mid |X_0^2u| \le 1\}}
  (D_1g_1)(u,X_0)\R(-X_1/X_0^2,X_1u) \dd u\] where $D_1g_1$ is the
  partial derivative of $g_1$ with respect to the first variable.
  Using the above bound for $\R(b_1,b_2)$, we obtain:
  \[\R_2(\ee) \ll \frac{X_2}{\Fthr}2^{\omega(\Cone)}.\]
  Summing this over all $\ee$ as in \eqref{eq:height_eta} while
  ignoring \eqref{eq:coprim_eta} which can only enlarge the sum, we
  obtain:
  \[\sum_{\ee} \R_2(\ee) \ll \sum_{\ee} \frac{X_2\cdot
    2^{\omega(\Cone)}}{\Fthr X_0^2} = \sum_{\ee} \frac{B\cdot
    2^{\omega(\Cone)}}{\sbase 2 1 2 2 2 2} \ll B\log B\]
  In the first step, we use $X_0 \le 1$.
\end{proof}

Let \[\Delta(n) = B^{-2/3}\!\!\!\sum_{\eta_i, \sbase 4 2 3 3 2 2 = n}
\frac{\vartheta(\ee) X_1X_2}{\Fthr} = \!\!\!\sum_{\eta_i,
  \sbase 4 2 3 3 2 2 = n}
\frac{\vartheta(\ee)(\sbase 4 2 3 3 2
  2)^{1/3}}{\sbase 1 1 1 1 1 1}.\]

In view of Lemma~\ref{lem:bijection}, the number of rational points of
bounded height on $\So$ can be estimated by summing the result of
Lemma~\ref{lem:sum_a1} over all suitable $\ee$. The error term is the
combination of the error terms in Lemmas~\ref{lem:sum_a2a3} and
\ref{lem:sum_a1}.

\begin{lemma}\label{lem:sum_eta}
  We have
  \[\NU = B^{2/3}\sum_{n \le B} \Delta(n)g_2((n/B)^{1/3}) + O(B (\log B)^3).\]
\end{lemma}

\section{Completion of the proof}\label{sec:proof}

We need an estimate for \[M(t) := \sum_{n \le t} \Delta(t).\] Consider
the Dirichlet series $F(s) := \sum_{n=1}^\infty \Delta(n)n^{-s}$.
Using \[F(s+1/3) = \sum_{\ee}
\frac{\vartheta(\ee)}{\eta_1^{4s+1}\eta_2^{2s+1}\eta_3^{3s+1}
  \eta_4^{3s+1}\eta_5^{2s+1}\eta_6^{2s+1}},\] we write $F(s+1/3) =
\prod_p F_p(s+1/3)$ as its Euler product. To obtain $F_p(s+1/3)$ for a
prime $p$, we need to restrict this sum to the terms in which all
$\eta_i$ are powers of $p$. Note that $\vartheta(\ee)$ is non-zero if
and only if the divisibility of $\eta_i$ by $p$ falls into one of the
twelve cases described in Table~\ref{tab:coprim_alpha}. The value of
$\vartheta(\ee)$ only depends on these cases.

Writing $F_p(s+1/3) = \sum_{i=1}^{11}F_{p,i}(s+1/3)$, we have for example:
\begin{align*}
  F_{p,0}(s+1/3) &= 1,\\
  F_{p,1}(s+1/3) &= \sum_{j=1}^\infty
  \frac{(1-1/p)(1-2/p)}{p^{j(4s+1)}}
  = \frac{(1-1/p)(1-2/p)}{p^{4s+1}-1},\\
  F_{p,7}(s+1/3) &= \sum_{j,k=1}^\infty
  \frac{(1-1/p)^2}{p^{j(4s+1)}p^{k(2s+1)}} =
  \frac{(1-1/p)^2}{(p^{4s+1}-1)(p^{2s+1}-1)}.
\end{align*}
The other cases are similiar, giving 
\begin{equation*}
  \begin{split}
    F_p(s+1/3) = &1+\frac{1-1/p}{p^{4s+1}-1}\bigg((1-2/p) +
    \frac{1-1/p}{p^{2s+1}-1} + 2\frac{1-1/p}{p^{3s+1}-1}\bigg)\\
    &+ \frac{1-1/p}{p^{2s+1}-1} + 2\frac{(1-1/p)^2}{p^{3s+1}-1} +
    2\frac{1-1/p}{p^{2s+1}-1} + 2\frac{(1-1/p)^2}{(p^{2s+1}-1)^2}.
  \end{split}
\end{equation*}

Defining \[E(s) := \zeta(4s+1)\zeta(3s+1)^2\zeta(2s+1)^3
\qquad\text{and}\qquad G(s) := F(s+1/3)/E(s),\] we see as in
\cite{math.NT/0504016} that the residue of $F(s)t^s/s$ at $s = 1/3$ is
\[\Res(t) = \frac{3 G(0)t^{1/3}Q_1(\log t)}{5! \cdot 4\cdot 2\cdot 3
  \cdot 3 \cdot 2 \cdot 2}\] for a monic $Q_1\in \mathbb R[x]$ of degree
$5$. By Lemma~\ref{lem:conjecture},
$
\alpha(S) = \frac{1}{5!\cdot 4\cdot 2\cdot 3 \cdot 3 \cdot 2
\cdot 2}.
$ 
 By a Tauberian argument
as in \cite[Lemma 13]{math.NT/0504016}:
\begin{lemma}
  $M(t) = \Res(t)+O(t^{1/3-\delta})$ for some $\delta > 0$.
\end{lemma}

By partial summation,
\[\sum_{n\le B}\Delta(n)g_2((n/B)^{1/3}) =\alpha(S)
G(0) B^{1/3}Q(\log B) \cdot 3\int_0^1 g_2(v) \dd v +
O(B^{\frac{1}{3}-\delta})\] for a monic polynomial $Q$ of degree 5.
We identify $\omega_H(S)$ from 
\[G(0) = \prod_p \left(1-\frac 1 p\right)^6 \left(1+\frac 6 p + \frac 1
  {p^2}\right),\,\, \text{ and } \omega_\infty = 3 \int_0^1 g_2(v) \dd v. \] 
Together with Lemma~\ref{lem:sum_eta}, this completes the proof of
Theorem~\ref{thm:main}.

\section{Equations of universal torsors}\label{sec:torsorexamples}

The simplest universal torsors are those which can be realized
as Zariski open subsets of the affine space. 
This happens if and only if the Del Pezzo surface is
toric. 

\begin{example}
  There are 20 types of singular Del Pezzo surfaces of degree $d \ge
  3$ whose universal torsor is an open subset of a hypersurface in
  $\AA^{13-d}$. For one example of each type\footnote{For the cubic
    $\Dfour$ case, the universal torsor of a different example is
    calculated in \cite[Section 4]{MR2029868}.}, the equation defining
  the universal torsor is listed in the following table. More details
  can be found in \cite{math.AG/0604194}.
 \[\begin{array}[ht]{|c|c|c|l|}
    \hline
    \text{degree} & \text{singularities} & \text{\# of lines} & 
    \text{defining equation}\\
    \hline\hline
    6 & \Aone & 3 & \eta_2\alpha_1+\eta_3\alpha_2+\eta_4\alpha_3 \\
    6 & \Atwo & 2 & \eta_2\alpha_1^2 +\eta_3\alpha_2+ \eta_4\alpha_3\\
    \hline
    5 & \Aone & 7 & \eta_2\eta_6+\eta_3\eta_7+\eta_4\eta_8 \\
    5 & \Atwo & 4 & 
    \eta_2\eta_5^2\eta_6 + \eta_3\alpha_1 + \eta_4\alpha_2 \\
    5 & \Athree & 2 &
    \eta_1\alpha_1^2 + \eta_3\eta_4^2\alpha_2 + \eta_5\alpha_3 \\
    5 & \Afour & 1 &
    \eta_1^2\eta_2\alpha_1^3 + \eta_4\alpha_2^2 + \eta_5\alpha_3 \\
    \hline
    4 & 3\Aone & 6 & 
    \eta_4\eta_5 + \eta_1\eta_6\eta_7 + \eta_8\eta_9 \\
    4 & \Atwo+\Aone & 6 & 
    \eta_5\eta_7 + \eta_1\eta_3\eta_9^2 + \eta_6\eta_8 \\
    4 & \Athree & 5 & 
    \eta_5\alpha + \eta_1\eta_4^2\eta_7 + \eta_3\eta_6^2\eta_8 \\
    4 & \Athree+\Aone & 3 & 
    \eta_6\alpha_2+\eta_7\alpha_1+\eta_1\eta_3\eta_4^2\eta_5^3\\
    4 & \Afour & 3 &
    \eta_5\alpha_1 + \eta_1\alpha_2^2 + \eta_3\eta_4^2\eta_6^3\eta_7 \\
    4 & \Dfour & 2 & 
    \eta_3\eta_5^2\alpha_2 + \eta_4\eta_6^2\alpha_3 + \eta_2\alpha_1^2\\
    4 & \Dfive & 1 & 
    \eta_3\alpha_1^2 + \eta_2\eta_6^2\alpha_3 + \eta_4\eta_5^2\alpha_2^3\\
    \hline
    3  & \Dfour & 6 & 
    \eta_2\eta_5^2\eta_8 + \eta_3\eta_6^2\eta_9 + \eta_4\eta_7^2\eta_{10}\\
    3 & \Athree + 2\Aone & 5 & 
    \eta_4\eta_6^2\eta_{10} + \eta_1\eta_2\eta_7^2 + \eta_8\eta_9 \\
    3 & 2\Atwo+\Aone & 5 & 
    \eta_3\eta_5\eta_7^2 + \eta_1\eta_6\eta_8 + \eta_9\eta_{10} \\
    3 & \Afour+\Aone & 4 & 
    \eta_1\eta_5\eta_8^2 + \eta_3\eta_4^2\eta_6^3\eta_9 + \eta_7\alpha \\
    3 & \Dfive & 3 & \eta_2\eta_6^2\alpha_2 + 
    \eta_4\eta_5^2\eta_7^3\eta_8 + \eta_3\alpha_1^2 \\
    3 & \Afive+\Aone & 2 & \eta_1^3\eta_2^2\eta_3\eta_7^4\eta_8 + 
    \eta_5\alpha_1^2 + \eta_6\alpha_2 \\
    3 & \Esix & 1 & \eta_4^2\eta_5\eta_7^3\alpha_3 + 
    \eta_2\alpha_2^2 + \eta_1^2\eta_3\alpha_1^3\\
    \hline
  \end{array}\]
\end{example}

\begin{example}[Cubic surface with $\Aone+\Athree$ singularities]
  This surface has 7 lines, 4 additional variables correspond to 
  exceptional curves of the desingularization. Its 9-dimensional
  universal torsor is a Zariski open subset of a complete intersection in 
  \[\AA^{11}=\Spec\ZZ[\eta_0, \dots, \eta_3, \mu_0, \dots, \mu_6]\]
  given by
  \[\eta_1\eta_2\mu_1\mu_2+\mu_4\mu_6+\mu_3\mu_5=0\qquad\text{and}\qquad
  \eta_0\eta_1\mu_2^2+\eta_3\mu_5\mu_6+\mu_0\mu_1=0.\]
  See \cite{cox_singular} for more details.
\end{example}

There are examples of universal torsors which are not
complete intersections, but have still been successfully used
in the context of Manin's conjecture:

\begin{example}[Cayley cubic]\label{ex:cayley}
  The Cayley cubic
  surface \[x_0x_1x_2+x_0x_1x_3+x_0x_2x_3+x_1x_2x_3=0\]
  (Figure~\ref{fig:cayley}) is a split singular cubic surface with
  four singularities $q_1, \dots, q_4$ of type $\Aone$ and nine lines.
  It is the blow-up of $\Ptwo$ in the 6 intersection points of 4 lines
  in general position.
  The universal torsor is an open subvariety of the variety in
  \[\AA^{13} =
  \Spec \ZZ[v_{12},v_{13},v_{14},y_1,y_2,y_3,y_4,z_{12},z_{13},z_{14},
  z_{23},z_{24},z_{25}]\] defined by six equations of the form
  \[z_{ik}z_{il}y_j+z_{jk}z_{jl}y_i = z_{ij}v_{ij}\]
  and three equations of the form
  \[v_{ij}v_{ik} = z_{il}^2y_jy_k-z_{jk}^2y_iy_l,\] where $\{i,j,k,l\}
  = \{1,2,3,4\}$ and \[z_{ij}=z_{ji}, \qquad v_{ij}=v_{ji},
  \qquad\text{and}\qquad v_{ij}=-v_{kl}.\] See \cite{cox_singular} for
  a proof. The variables $y_i$ correspond to the four exceptional
  curves $E_i$ obtained by blowing up $q_i$, $z_{ij}$ correspond to
  the six lines $m_{ij}$ through two of the singularities, and
  $v_{ij}$ correspond to the other three lines $\ell_{ij}$. The first
  six equations can be interpreted in connection with the projection
  from $m_{ij}$, and the other three equations are connected to the
  projection from $\ell_{ij}$.
  
  Upper and lower bounds of the expected order of magnitude have been
  established in \cite{MR2075628}.
\end{example}

\begin{example}[Smooth degree 5 Del Pezzo surface]
  The blow-up of $\Ptwo$ in \[(1:0:0), \qquad (0:1:0), \qquad 
  (0:0:1), \qquad (1:1:1)\] is a split smooth Del Pezzo surface of degree 5.
  Its universal torsor is an open subset of the variety defined by the
  following five equations in ten variables:
  \begin{align*}
    \lambda_{12}\eta_2 - \lambda_{13}\eta_3 + \lambda_{14}\eta_4 &= 0\\
    \lambda_{12}\eta_1 - \lambda_{23}\eta_3 + \lambda_{24}\eta_4 &= 0\\
    \lambda_{13}\eta_1 - \lambda_{23}\eta_2 + \lambda_{34}\eta_4 &= 0\\
    \lambda_{14}\eta_1 - \lambda_{24}\eta_2 + \lambda_{34}\eta_3 &= 0\\
    \lambda_{12}\lambda_{34} - \lambda_{13}\lambda_{24} +
    \lambda_{14}\lambda_{23} &= 0
  \end{align*}
  The asymptotic formula (\ref{eqn:asym}) has been established in
  \cite{MR1909606}.
\end{example}

To illustrate some of the difficulties in proving Conjecture
\ref{conj:manin} for a smooth split cubic surface, we now write down
equations for its universal torsor (up to radical).
  
\begin{example}[Smooth cubic surfaces]
  Let $S$ be the blow-up of $\Ptwo$ in \[(1:0:0), \,\, (0:1:0), \,\,
  (0:0:1),\, \, (1:1:1), \,\, (1:a:b),\, \, (1:c:d),\] in general
  position.  Conjecturally, the universal torsor is an open subset of
  the intersection of 81 quadrics in 27-dimensional space $\Spec
  \ZZ[\eta_i, \mu_{i,j}, \lambda_i]$, where
  \begin{itemize}
  \item $\eta_1, \dots, \eta_6$ correspond to the preimages of the points,
  \item $\mu_{i,j}$ ($i < j \in \{1,\dots, 6\}$) correspond to the 15 lines
    $m_{i,j}$ through two of the points,
  \item $\lambda_1, \dots, \lambda_6$ correspond to the conics $Q_i$
    through five of the six points,
  \end{itemize}
  and relations arise from conic bundle structures on $S$. Batyrev and Popov
  proved that the above variables are indeed generators and that the relations
  give the universal torsor, up to radical \cite{MR2029863}. 

  We now write down these equations explicitly (see
  \cite{math.AG/0603111} for more details).  The 81 defining quadrics
  occur in sets of three. These 27 triples correspond to projections
  from the 27 lines on $S$. We use
  \[E:=(b-1)(c-1)-(a-1)(d-1)\,\, \text{ and }\,\, F:=bc-ad\]
  to simplify the equations.
\begin{align*}
\qa{Q_1}&=-\eta_2\mu_{1,2}-\eta_3\mu_{1,3}+\eta_4\mu_{1,4}\\
\qb{Q_1}&=-a\eta_2\mu_{1,2}-b\eta_3\mu_{1,3}+\eta_5\mu_{1,5}\\
\qb{Q_1}&=-c\eta_2\mu_{1,2}-d\eta_3\mu_{1,3}+\eta_6\mu_{1,6}
\displaybreak[0]\\[\baselineskip]
\qa{Q_2}&=\eta_1\mu_{1,2}-\eta_3\mu_{2,3}+\eta_4\mu_{2,4}\\
\qb{Q_2}&=\eta_1\mu_{1,2}-b\eta_3\mu_{2,3}+\eta_5\mu_{2,5}\\
\qc{Q_2}&=\eta_1\mu_{1,2}-d\eta_3\mu_{2,3}+\eta_6\mu_{2,6}
\displaybreak[0]\\[\baselineskip]
\qa{Q_3}&=\eta_1\mu_{1,3}+\eta_2\mu_{2,3}+\eta_4\mu_{3,4}\\
\qb{Q_3}&=\eta_1\mu_{1,3}+a\eta_2\mu_{2,3}+\eta_5\mu_{3,5}\\
\qc{Q_3}&=\eta_1\mu_{1,3}+c\eta_2\mu_{2,3}+\eta_6\mu_{3,6}
\displaybreak[0]\\[\baselineskip]
\qa{Q_4}&=\eta_1\mu_{1,4}+\eta_2\mu_{2,4}+\eta_3\mu_{3,4}\\
\qb{Q_4}&=(1-b)\eta_1\mu_{1,4}+(a-b)\eta_2\mu_{2,4}+\eta_5\mu_{4,5}\\
\qc{Q_4}&=(1-d)\eta_1\mu_{1,4}+(c-d)\eta_2\mu_{2,4}+\eta_6\mu_{4,6}
\displaybreak[0]\\[\baselineskip]
\qa{Q_5}&=1/b\eta_1\mu_{1,5}+a/b\eta_2\mu_{2,5}+\eta_3\mu_{3,5}\\
\qb{Q_5}&=(1-b)/b\eta_1\mu_{1,5}+(a-b)/b\eta_2\mu_{2,5}+\eta_4\mu_{4,5}\\
\qc{Q_5}&=(b-d)/b\eta_1\mu_{1,5}+F/b\eta_2\mu_{2,5}+\eta_6\mu_{5,6}
\displaybreak[0]\\[\baselineskip]
\qa{Q_6}&=1/d\eta_1\mu_{1,6}+c/d\eta_2\mu_{2,6}+\eta_3\mu_{3,6}\\
\qb{Q_6}&=(1-d)/d\eta_1\mu_{1,6}+(c-d)/d\eta_2\mu_{2,6}+\eta_4\mu_{4,6}\\
\qc{Q_6}&=(b-d)/d\eta_1\mu_{1,6}+F/d\eta_2\mu_{2,6}+\eta_5\mu_{5,6}
\displaybreak[0]\\[\baselineskip]
\qa{m_{1,2}}&=\mu_{4,5}\mu_{3,6}-\mu_{3,5}\mu_{4,6}+\mu_{3,4}\mu_{5,6}\\
\qb{m_{1,2}}&=(b-d)\mu_{3,5}\mu_{4,6}+(d-1)\mu_{3,4}\mu_{5,6}+\eta_2\lambda_1\\
\qc{m_{1,2}}&=F\mu_{3,5}\mu_{4,6}+a(d-c)\mu_{3,4}\mu_{5,6}+\eta_1\lambda_2
\displaybreak[0]\\[\baselineskip]
\qa{m_{1,3}}&=\mu_{4,5}\mu_{2,6}-\mu_{2,5}\mu_{4,6}+\mu_{2,4}\mu_{5,6}\\
\qb{m_{1,3}}&=(c-a)\mu_{2,5}\mu_{4,6}+(1-c)\mu_{2,4}\mu_{5,6}+\eta_3\lambda_1\\
\qc{m_{1,3}}&=-F\mu_{2,5}\mu_{4,6}+b(c-d)\mu_{2,4}\mu_{5,6}+\eta_1\lambda_3
\displaybreak[0]\\[\baselineskip]
\qa{m_{2,3}}&=\mu_{4,5}\mu_{1,6}-\mu_{1,5}\mu_{4,6}+\mu_{1,4}\mu_{5,6}\\
\qb{m_{2,3}}&=(a-c)\mu_{1,5}\mu_{4,6}+a(c-1)\mu_{1,4}\mu_{5,6}+\eta_3\lambda_2\\
\qc{m_{2,3}}&=(b-d)\mu_{1,5}\mu_{4,6}+b(d-1)\mu_{1,4}\mu_{5,6}+\eta_2\lambda_3
\displaybreak[0]\\[\baselineskip]
\qa{m_{1,4}}&=\mu_{3,5}\mu_{2,6}-\mu_{2,5}\mu_{3,6}+\mu_{2,3}\mu_{5,6}\\
\qb{m_{1,4}}&=-E\mu_{2,5}\mu_{3,6}+(b-1)(c-1)\mu_{2,3}\mu_{5,6}+\eta_4\lambda_1\\
\qc{m_{1,4}}&=-F\mu_{2,5}\mu_{3,6}+bc\mu_{2,3}\mu_{5,6}+\eta_1\lambda_4
\displaybreak[0]\\[\baselineskip]
\qa{m_{2,4}}&=\mu_{3,5}\mu_{1,6}-\mu_{1,5}\mu_{3,6}+\mu_{1,3}\mu_{5,6}\\
\qb{m_{2,4}}&=E\mu_{1,5}\mu_{3,6}+(a-b)(c-1)\mu_{1,3}\mu_{5,6}+\eta_4\lambda_2\\
\qc{m_{2,4}}&=(b-d)\mu_{1,5}\mu_{3,6}-b\mu_{1,3}\mu_{5,6}+\eta_2\lambda_4
\displaybreak[0]\\[\baselineskip]
\qa{m_{3,4}}&=\mu_{2,5}\mu_{1,6}-\mu_{1,5}\mu_{2,6}+\mu_{1,2}\mu_{5,6}\\
\qb{m_{3,4}}&=-E\mu_{1,5}\mu_{2,6}+(a-b)(1-d)\mu_{1,2}\mu_{5,6}+\eta_4\lambda_3\\
\qc{m_{3,4}}&=(c-a)\mu_{1,5}\mu_{2,6}+a\mu_{1,2}\mu_{5,6}+\eta_3\lambda_4
\displaybreak[0]\\[\baselineskip]
\qa{m_{1,5}}&=\mu_{3,4}\mu_{2,6}-\mu_{2,4}\mu_{3,6}+\mu_{2,3}\mu_{4,6}\\
\qb{m_{1,5}}&=-E\mu_{2,4}\mu_{3,6}+(a-c)(1-b)\mu_{2,3}\mu_{4,6}+\eta_5\lambda_1\\
\qc{m_{1,5}}&=(d-c)\mu_{2,4}\mu_{3,6}+c\mu_{2,3}\mu_{4,6}+\eta_1\lambda_5
\displaybreak[0]\\[\baselineskip]
\qa{m_{2,5}}&=\mu_{3,4}\mu_{1,6}-\mu_{1,4}\mu_{3,6}+\mu_{1,3}\mu_{4,6}\\
\qb{m_{2,5}}&=aE\mu_{1,4}\mu_{3,6}+(a-b)(c-a)\mu_{1,3}\mu_{4,6}+\eta_5\lambda_2\\
\qc{m_{2,5}}&=(1-d)\mu_{1,4}\mu_{3,6}-\mu_{1,3}\mu_{4,6}+\eta_2\lambda_5
\displaybreak[0]\\[\baselineskip]
\qa{m_{3,5}}&=\mu_{2,4}\mu_{1,6}-\mu_{1,4}\mu_{2,6}+\mu_{1,2}\mu_{4,6}\\
\qb{m_{3,5}}&=-bE\mu_{1,4}\mu_{2,6}+(a-b)(b-d)\mu_{1,2}\mu_{4,6}+\eta_5\lambda_3\\
\qc{m_{3,5}}&=(c-1)\mu_{1,4}\mu_{2,6}+\mu_{1,2}\mu_{4,6}+\eta_3\lambda_5
\displaybreak[0]\\[\baselineskip]
\qa{m_{4,5}}&=\mu_{2,3}\mu_{1,6}-\mu_{1,3}\mu_{2,6}+\mu_{1,2}\mu_{3,6}\\
\qb{m_{4,5}}&=b(c-a)\mu_{1,3}\mu_{2,6}+a(b-d)\mu_{1,2}\mu_{3,6}+\eta_5\lambda_4\\
\qc{m_{4,5}}&=(c-1)\mu_{1,3}\mu_{2,6}+(1-d)\mu_{1,2}\mu_{3,6}+\eta_4\lambda_5
\displaybreak[0]\\[\baselineskip]
\qa{m_{1,6}}&=\mu_{3,4}\mu_{2,5}-\mu_{2,4}\mu_{3,5}+\mu_{2,3}\mu_{4,5}\\
\qb{m_{1,6}}&=-E\mu_{2,4}\mu_{3,5}+(a-c)(1-d)\mu_{2,3}\mu_{4,5}+\eta_6\lambda_1\\
\qc{m_{1,6}}&=(b-a)\mu_{2,4}\mu_{3,5}+a\mu_{2,3}\mu_{4,5}+\eta_1\lambda_6
\displaybreak[0]\\[\baselineskip]
\qa{m_{2,6}}&=\mu_{3,4}\mu_{1,5}-\mu_{1,4}\mu_{3,5}+\mu_{1,3}\mu_{4,5}\\
\qb{m_{2,6}}&=cE\mu_{1,4}\mu_{3,5}+(a-c)(d-c)\mu_{1,3}\mu_{4,5}+\eta_6\lambda_2\\
\qc{m_{2,6}}&=(1-b)\mu_{1,4}\mu_{3,5}-\mu_{1,3}\mu_{4,5}+\eta_2\lambda_6
\displaybreak[0]\\[\baselineskip]
\qa{m_{3,6}}&=\mu_{2,4}\mu_{1,5}-\mu_{1,4}\mu_{2,5}+\mu_{1,2}\mu_{4,5}\\
\qb{m_{3,6}}&=-dE\mu_{1,4}\mu_{2,5}+(d-b)(d-c)\mu_{1,2}\mu_{4,5}+\eta_6\lambda_3\\
\qc{m_{3,6}}&=(a-1)\mu_{1,4}\mu_{2,5}+\mu_{1,2}\mu_{4,5}+\eta_3\lambda_6
\displaybreak[0]\\[\baselineskip]
\qa{m_{4,6}}&=\mu_{2,3}\mu_{1,5}-\mu_{1,3}\mu_{2,5}+\mu_{1,2}\mu_{3,5}\\
\qb{m_{4,6}}&=d(c-a)\mu_{1,3}\mu_{2,5}+c(b-d)\mu_{1,2}\mu_{3,5}+\eta_6\lambda_4\\
\qc{m_{4,6}}&=(a-1)\mu_{1,3}\mu_{2,5}+(1-b)\mu_{1,2}\mu_{3,5}+\eta_4\lambda_6
\displaybreak[0]\\[\baselineskip]
\qa{m_{5,6}}&=\mu_{2,3}\mu_{1,4}-\mu_{1,3}\mu_{2,4}+\mu_{1,2}\mu_{3,4}\\
\qb{m_{5,6}}&=d(c-1)\mu_{1,3}\mu_{2,4}+c(1-d)\mu_{1,2}\mu_{3,4}+\eta_6\lambda_5\\
\qc{m_{5,6}}&=b(a-1)\mu_{1,3}\mu_{2,4}+a(1-b)\mu_{1,2}\mu_{3,4}+\eta_5\lambda_6
\displaybreak[0]\\[\baselineskip]
\qa{E_1}&=(d-b)/E\mu_{1,2}\lambda_2+(c-a)/E\mu_{1,3}\lambda_3+\mu_{1,4}\lambda_4\\
\qb{E_1}&=(d-1)/E\mu_{1,2}\lambda_2+(c-1)/E\mu_{1,3}\lambda_3+\mu_{1,5}\lambda_5\\
\qc{E_1}&=(b-1)/E\mu_{1,2}\lambda_2+(a-1)/E\mu_{1,3}\lambda_3+\mu_{1,6}\lambda_6
\displaybreak[0]\\[\baselineskip]
\qa{E_2}&=F/E\mu_{1,2}\lambda_1+(c-a)/E\mu_{2,3}\lambda_3+\mu_{2,4}\lambda_4\\
\qb{E_2}&=(c-d)/E\mu_{1,2}\lambda_1+(c-1)/E\mu_{2,3}\lambda_3+\mu_{2,5}\lambda_5\\
\qc{E_2}&=(a-b)/E\mu_{1,2}\lambda_1+(a-1)/E\mu_{2,3}\lambda_3+\mu_{2,6}\lambda_6
\displaybreak[0]\\[\baselineskip]
\qa{E_3}&=F/E\mu_{1,3}\lambda_1+(b-d)/E\mu_{2,3}\lambda_2+\mu_{3,4}\lambda_4\\
\qb{E_3}&=(c-d)/E\mu_{1,3}\lambda_1+(1-d)/E\mu_{2,3}\lambda_2+\mu_{3,5}\lambda_5\\
\qc{E_3}&=(a-b)/E\mu_{1,3}\lambda_1+(1-b)/E\mu_{2,3}\lambda_2+\mu_{3,6}\lambda_6
\displaybreak[0]\\[\baselineskip]
\qa{E_4}&=F/(a-c)\mu_{1,4}\lambda_1+(b-d)/(a-c)\mu_{2,4}\lambda_2+\mu_{3,4}\lambda_3\\
\qb{E_4}&=c/(a-c)\mu_{1,4}\lambda_1+1/(a-c)\mu_{2,4}\lambda_2+\mu_{4,5}\lambda_5\\
\qc{E_4}&=a/(a-c)\mu_{1,4}\lambda_1+1/(a-c)\mu_{2,4}\lambda_2+\mu_{4,6}\lambda_6
\displaybreak[0]\\[\baselineskip]
\qa{E_5}&=(d-c)/(c-1)\mu_{1,5}\lambda_1+(d-1)/(c-1)\mu_{2,5}\lambda_2+\mu_{3,5}\lambda_3\\
\qb{E_5}&=-c/(c-1)\mu_{1,5}\lambda_1-1/(c-1)\mu_{2,5}\lambda_2+\mu_{4,5}\lambda_4\\
\qc{E_5}&=-1/(c-1)\mu_{1,5}\lambda_1-1/(c-1)\mu_{2,5}\lambda_2+\mu_{5,6}\lambda_6
\displaybreak[0]\\[\baselineskip]
\qa{E_6}&=(b-a)/(a-1)\mu_{1,6}\lambda_1+(b-1)/(a-1)\mu_{2,6}\lambda_2+\mu_{3,6}\lambda_3\\
\qb{E_6}&=-a/(a-1)\mu_{1,6}\lambda_1-1/(a-1)\mu_{2,6}\lambda_2+\mu_{4,6}\lambda_4\\
\qc{E_6}&=-1/(a-1)\mu_{1,6}\lambda_1-1/(a-1)\mu_{2,6}\lambda_2+\mu_{5,6}\lambda_5
\end{align*}
\end{example}

\

In general, the dimension $k$ of the ambient space $\AA^k$ of the
universal torsor is at least as large as the number of lines on the
surface plus the number of exceptional curves of its
desingularization, while the dimension of the universal torsor only
depends on the degree of the surface, so that the number of equations
must grow with $k$.

Heuristically, the complexity of universal torsors should be dictated
by the following considerations:
\begin{itemize}
\item The dimension of the universal torsor of split Del Pezzo
  surfaces $S$ is $12-d$, where $d$ is the degree of $S$.
\item For smooth Del Pezzo surfaces, the number of lines is bigger in
  smaller degrees (e.g., 10 lines in degree 5, and 27 lines in degree
  3).
\item Singular surfaces have less lines than smooth surfaces.
\item The number of lines is higher in cases with ``few mild''
  singularities (e.g., for cubics: $\Aone$ with 21 lines, $\Atwo$ with
  15 lines), while it is low for ``bad'' singularities (e.g., 1 for
  the $\Esix$ cubic, 2 for the $\Afive+\Aone$ cubic).
\end{itemize}
Therefore, we expect universal torsors over surfaces which have
low degree, are smooth or have mild singularities to be more complex
than torsors over surfaces in large degree, or with complicated
singularities.

\bibliographystyle{alpha}

\bibliography{qd4}

\begin{thebibliography}{FMT89}

\bibitem[BB]{math.NT/0412086}
{R. de la} Bret{\`e}che and T.~D. Browning.
\newblock {On Manin's conjecture for singular del Pezzo surfaces of degree
  four, I}.
\newblock {\em Mich. Math. J.}, to appear.

\bibitem[BB05]{math.NT/0502510}
{R. de la} Bret{\`e}che and T.~D. Browning.
\newblock {On Manin's conjecture for singular del Pezzo surfaces of degree
  four, II}, arXiv:math.NT/0502510, 2005.

\bibitem[BBD05]{math.NT/0509370}
{R. de la} Bret{\`e}che, T.~D. Browning, and U.~Derenthal.
\newblock {On Manin's conjecture for a certain singular cubic surface},
  arXiv:math.NT/0509370, 2005.

\bibitem[BF04]{MR2099200}
{R. de la} Bret{\`e}che and {\'E}.~Fouvry.
\newblock L'\'eclat\'e du plan projectif en quatre points dont deux
  conjugu\'es.
\newblock {\em J. Reine Angew. Math.}, 576:63--122, 2004.

\bibitem[Bir62]{MR0150129}
B.~J. Birch.
\newblock Forms in many variables.
\newblock {\em Proc. Roy. Soc. Ser. A}, 265:245--263, 1961/1962.

\bibitem[BM90]{MR1032922}
V.~V. Batyrev and Y.~I. Manin.
\newblock Sur le nombre des points rationnels de hauteur born\'e des
  vari\'et\'es alg\'ebriques.
\newblock {\em Math. Ann.}, 286(1-3):27--43, 1990.

\bibitem[BP04]{MR2029863}
V.~V. Batyrev and O.~N. Popov.
\newblock The {C}ox ring of a del {P}ezzo surface.
\newblock In {\em Arithmetic of higher-dimensional algebraic varieties (Palo
  Alto, CA, 2002)}, volume 226 of {\em Progr. Math.}, pages 85--103.
  Birkh\"auser Boston, Boston, MA, 2004.

\bibitem[Bre98]{MR2000b:11074}
R.~de~la Bret{\`e}che.
\newblock Sur le nombre de points de hauteur born\'ee d'une certaine surface
  cubique singuli\`ere.
\newblock {\em Ast\'erisque}, (251):51--77, 1998.
\newblock Nombre et r\'epartition de points de hauteur born\'ee (Paris, 1996).

\bibitem[Bre02]{MR1909606}
R.~de~la Bret{\`e}che.
\newblock Nombre de points de hauteur born\'ee sur les surfaces de del {P}ezzo
  de degr\'e 5.
\newblock {\em Duke Math. J.}, 113(3):421--464, 2002.

\bibitem[Bro04]{math.NT/0404245}
T.~D. Browning.
\newblock {The density of rational points on a certain singular cubic surface},
  arXiv:math.NT/0404245, 2004.

\bibitem[Bro05]{math.NT/0511041}
T.~D. Browning.
\newblock {An overview of Manin's conjecture for del Pezzo surfaces},
  arXiv:math.NT/0511041, 2005.

\bibitem[BT98]{MR1620682}
V.~V. Batyrev and Y.~Tschinkel.
\newblock Manin's conjecture for toric varieties.
\newblock {\em J. Algebraic Geom.}, 7(1):15--53, 1998.

\bibitem[CLT02]{MR1906155}
A.~Chambert-Loir and Y.~Tschinkel.
\newblock On the distribution of points of bounded height on equivariant
  compactifications of vector groups.
\newblock {\em Invent. Math.}, 148(2):421--452, 2002.

\bibitem[CTS87]{MR89f:11082}
J.-L. Colliot-Th{\'e}l{\`e}ne and J.-J. Sansuc.
\newblock La descente sur les vari\'et\'es rationnelles. {II}.
\newblock {\em Duke Math. J.}, 54(2):375--492, 1987.

\bibitem[Der05]{math.NT/0504016}
U.~Derenthal.
\newblock {Manin's conjecture for a certain singular cubic surface},
  arXiv:math.NT/0504016, 2005.

\bibitem[Der06a]{del_pezzo_alpha}
U.~Derenthal.
\newblock {On a constant arising in Manin's conjecture for Del Pezzo surfaces},
  2006.

\bibitem[Der06b]{cox_singular}
U.~Derenthal.
\newblock {On the Cox ring of singular Del Pezzo surfaces}, 2006.

\bibitem[Der06c]{math.AG/0603111}
U.~Derenthal.
\newblock {On the Cox ring of Del Pezzo surfaces}, arXiv:math.AG/0603111, 2006.

\bibitem[Der06d]{math.AG/0604194}
U.~Derenthal.
\newblock {Singular Del Pezzo surfaces whose universal torsors are
  hypersurfaces}, arXiv:math.AG/0604194, 2006.

\bibitem[Dol03]{MR2004511}
I.~Dolgachev.
\newblock {\em Lectures on invariant theory}, volume 296 of {\em London
  Mathematical Society Lecture Note Series}.
\newblock Cambridge University Press, Cambridge, 2003.

\bibitem[FMT89]{MR89m:11060}
J.~Franke, Y.~I. Manin, and Y.~Tschinkel.
\newblock Rational points of bounded height on {F}ano varieties.
\newblock {\em Invent. Math.}, 95(2):421--435, 1989.

\bibitem[Har77]{MR0463157}
R.~Hartshorne.
\newblock {\em Algebraic geometry}.
\newblock Springer-Verlag, New York, 1977.
\newblock Graduate Texts in Mathematics, No. 52.

\bibitem[HB03]{MR2075628}
D.~R. Heath-Brown.
\newblock The density of rational points on {C}ayley's cubic surface.
\newblock In {\em Proceedings of the Session in Analytic Number Theory and
  Diophantine Equations}, volume 360 of {\em Bonner Math. Schriften}, page~33,
  Bonn, 2003. Univ. Bonn.

\bibitem[HK00]{MR2001i:14059}
Y.~Hu and S.~Keel.
\newblock Mori dream spaces and {GIT}.
\newblock {\em Michigan Math. J.}, 48:331--348, 2000.
\newblock Dedicated to William Fulton on the occasion of his 60th birthday.

\bibitem[HT04]{MR2029868}
B.~Hassett and Y.~Tschinkel.
\newblock Universal torsors and {C}ox rings.
\newblock In {\em Arithmetic of higher-dimensional algebraic varieties (Palo
  Alto, CA, 2002)}, volume 226 of {\em Progr. Math.}, pages 149--173.
  Birkh\"auser Boston, Boston, MA, 2004.

\bibitem[Man86]{MR833513}
Yu.~I. Manin.
\newblock {\em Cubic forms}, volume~4 of {\em North-Holland Mathematical
  Library}.
\newblock North-Holland Publishing Co., Amsterdam, second edition, 1986.
\newblock Algebra, geometry, arithmetic, Translated from the Russian by M.
  Hazewinkel.

\bibitem[Pey95]{MR1340296}
E.~Peyre.
\newblock Hauteurs et mesures de {T}amagawa sur les vari\'et\'es de {F}ano.
\newblock {\em Duke Math. J.}, 79(1):101--218, 1995.

\bibitem[Pey98]{MR1679842}
E.~Peyre.
\newblock Terme principal de la fonction z\^eta des hauteurs et torseurs
  universels.
\newblock {\em Ast\'erisque}, (251):259--298, 1998.
\newblock Nombre et r\'epartition de points de hauteur born\'ee (Paris, 1996).

\bibitem[Sal98]{MR1679841}
P.~Salberger.
\newblock Tamagawa measures on universal torsors and points of bounded height
  on {F}ano varieties.
\newblock {\em Ast\'erisque}, (251):91--258, 1998.
\newblock Nombre et r\'epartition de points de hauteur born\'ee (Paris, 1996).

\end{thebibliography}

\end{document}